\documentclass[a4paper,12pt]{amsart} 
\usepackage{amsfonts, amssymb, amsmath}
\usepackage{longtable}
\usepackage[all]{xy}
\usepackage{epsfig}

\theoremstyle{plain}
\newtoks\thehProclaim
\newtheorem*{Proclaim}{\the\thehProclaim}
\newenvironment{proclaim}[1]{\thehProclaim{#1}\begin{Proclaim}}{\end{Proclaim}}

\def\Hom{\text{\rm Hom}}

\def\rk{\operatorname{rk}}

\def\e{\varepsilon}
\def\a{\alpha}
\def\b{\beta}

\def\A{\operatorname{A}}
\def\B{\operatorname{B}}
\def\C{\operatorname{C}}
\def\D{\operatorname{D}}

\def\K{\operatorname{K}}

\def\SL{\operatorname{SL}}
\def\SD{\operatorname{SD}}
\def\GL{\operatorname{GL}}

\def\St{\operatorname{St}}

\def\ad{\operatorname{ad}}
\def\sr{\operatorname{sr}}
\def\asr{\operatorname{asr}}

\def\sic{\operatorname{sc}}

\def\Int{{\mathbb Z}}

\def\Co{{\mathbb C}}

\begin{document}

\title[Gauss decomposition for Chevalley groups]
{Gauss decomposition for\\ Chevalley groups, revisited}
\author{A.~Smolensky, B.~Sury, N.~Vavilov}

\address{\newline
St. Petersburg State University\newline
University pr. 28,
Peterhof, \newline 198504 St. Petersburg, Russia
\bigskip
\newline
Indian Statistical Institute\newline 8th Mile Mysore Road\newline
Bangalore 560059, India\newline\bigskip}

\email {
\newline
andrei.smolensky@gmail.com\newline
surybang@gmail.com\newline
nikolai-vavilov@yandex.ru}
\date{20 July 2011}
\keywords{ Chevalley groups, elementary Chevalley groups, triangular
factorisations, rings of stable rank 1, parabolic subgroups, Gauss
decomposition, commutator width}

\thanks
{The work of the first and the third author
was supported by the State Financed task project 6.38.74.2011
at the Saint Petersburg State University
``Structure theory and geometry of algebraic groups, and their
applications in representation theory and algebraic $\K$-theory''.
The work of the second and the third author was conducted
in the framework of a joint Russian--Indian project RFFI 10-01-92651
``Higher composition laws, algebraic $\K$-theory and exceptional
groups'' (SPbGU).
The second author is grateful to the Saint Petersburg Department of
Steklov Mathematical Institute and to the Saint Petersburg State
University for the invitation to visit Saint Petersburg in
May--June 2011.
The research of the third author was carried out in the
framework of the RFFI projects 09-01-00784, ``Efficient
generation in groups of Lie type'' (PDMI RAS),
09-01-00878 ``Overgroups of reductive groups in algebraic
groups over rings'' (SPbGU), 10-01-90016
``The study of structure of forms of reductive groups and behaviour
of small unipotent elements in representations of algebraic groups''
(SPbGU). Apart from that, he was partially supported by the
RFFI projects 09-01-00762 (Siberian Federal University),
09-01-91333 (POMI RAS) and 11-01-00756 (RGPU)}

\maketitle

\begin{abstract}
In the 1960's Noboru Iwahori and Hideya Matsumoto, Eiichi
Abe and
Kazuo Suzuki, and Michael Stein discovered that Chevalley groups
$G=G(\Phi,R)$ over a semilocal ring admit remarkable Gauss
decomposition $G=TUU^-U$, where $T=T(\Phi,R)$ is a split maximal
torus, whereas $U=U(\Phi,R)$ and $U^-=U^-(\Phi,R)$ are unipotent
radicals of two opposite Borel subgroups $B=B(\Phi,R)$ and
$B^-=B^-(\Phi,R)$ containing $T$. It follows from the classical work
of Hyman Bass and Michael Stein that for classical groups Gauss
decomposition holds under weaker assumptions such as $\sr(R)=1$ or
$\asr(R)=1$. Later the second author noticed that condition
$\sr(R)=1$ is necessary for Gauss decomposition. Here, we show that
a slight variation of Tavgen's rank reduction theorem implies that
for the elementary group $E(\Phi,R)$ condition $\sr(R)=1$ is also
sufficient for Gauss decomposition. In other words, $E=HUU^-U$,
where $H=H(\Phi,R)=T\cap E$. This surprising result shows that
stronger conditions on the ground ring, such as being semi-local,
$\asr(R)=1$, $\sr(R,\Lambda)=1$, etc., were only needed to guarantee
that for simply connected groups $G=E$, rather than to verify the
Gauss decomposition itself.
\end{abstract}

\goodbreak

Let $\Phi$ be a reduced irreducible root system, $R$
be a commutative ring with 1 and $G(\Phi,R)$ be a Chevalley
group of type $\Phi$ over $R$. We fix a split maximal
torus $T(\Phi,R)$ in $G(\Phi,R)$ and a pair $B(\Phi,R)$
and $B^-(\Phi,R)$ of opposite Borel subgroups containing
$T(\Phi,R)$. Further, let $U(\Phi,R)$ and $U^-(\Phi,R)$
be the unipotent radicals of $B(\Phi,R)$ and $B^-(\Phi,R)$,
respectively.
\par
The whole theory of Chevalley groups over semi-local rings
rests upon the following analogue of Gauss decomposition
established by Noboru Iwahori and Hideya Matsumoto \cite{IM},
Eiichi Abe and Kazuo Suzuki \cite{abe,AS}, and
by Michael Stein \cite{stein73}. In fact, it plays the same role
in this case, as Bruhat decomposition does over fields.
Let\/ $R$ be a semi-local ring. Then one has the following
decomposition
$$ G(\Phi,R)=T(\Phi,R)U(\Phi,R)U^-(\Phi,R)U(\Phi,R). $$
\par
For the simply connected Chevalley group $G_{\sic}(\Phi,R)$
its maximal torus $T_{\sic}(\Phi,R)$ is contained in the
elementary subgroup
$$ E(\Phi,R)=\big\langle U(\Phi,R),U^-(\Phi,R)\big\rangle, $$
\noindent
generated by $U(\Phi,R)$ and $U^-(\Phi,R)$. In particular,
Gauss decomposition implies that for a simply connected
group over a semi-local ring one has
$G_{\sic}(\Phi,R)=E_{\sic}(\Phi,R)$. In other words,
for semi-local rings
$$ K_1(\Phi,R)=G_{\sic}(\Phi,R)/E_{\sic}(\Phi,R) $$
\noindent
is trivial.
\par
In general, when the group is not simply connected or the
ring $R$ is not semi-local, elementary subgroup $E(\Phi,R)$
can be strictly smaller than the Chevalley group $G(\Phi,R)$
itself. In fact, for a non simply connected group even the
subgroup $H(\Phi,R)$ spanned by semi-simple root elements
$h_{\a}(\e)$, $\a\in\Phi$, $\e\in R^*$, where $R^*$ is
the multiplicative group of the ring $R$. Clearly,
$$ H(\Phi,R)=T(\Phi,R)\cap E(\Phi,R), $$
\noindent
can be strictly smaller than the torus $T(\Phi,R)$ itself.
\par
In \cite{vavilov84} the third author observed that condition $\sr(R)=1$
is {\it necessary\/} for Gauss decomposition to hold for
a Chevalley group $G(\Phi,R)$ over a ring $R$ and made the
following remark:
``One might hope that the condition $\sr(R)=1$ is also
{\it sufficient\/} to prove that Chevalley groups of all types
over $R$ there admit Gauss decomposition (some
experts believe that this is rather unlikely).''
\par
In the present paper, which is a sequel to our paper \cite{VSS}, we
show that a slight modification of the same argument by Oleg Tavgen
\cite{tavgen90}, immediately gives the following surprising result,
asserting that condition $\sr(R)=1$ is necessary and {\it
sufficient\/} for the {\it elementary\/} Chevalley group $E(\Phi,R)$
to admit Gauss decompostion.
\begin{proclaim}
{Theorem 1} Let\/ $\Phi$ be a reduced irreducible root system
and\/ $R$ be a commutative ring such that $\sr(R)=1$. Then
the elementary Chevalley group\/ $E(\Phi,R)$ admits Gauss
decomposition
$$ E(\Phi,R)=H(\Phi,R)U(\Phi,R)U^-(\Phi,R)U(\Phi,R). $$
\noindent
Conversely, if Gauss decomposition holds for some\/
{\rm[}elementary\/{\rm]} Chevalley group, then\/ $\sr(R)=1$.
\end{proclaim}
The proof of this result follows {\it exactly\/} the same lines as
the proof of Theorem~1 in \cite{VSS}\footnote{The statement of
Theorem 1 in the Russian original of \cite{VSS} contains a very
unfortunate misprint. The formula there reads as the
unitriangular factorisation of length 4 for the simply connected
Chevalley group $G(\Phi,R)$. Of course, the rest of that paper
and the proofs there discuss such a factorisation for the
{\it elementary\/} Chevalley group $E(\Phi,R)$. This inconsistency
is corrected in the English version.}.
Now, we are interested not in
unitriangular factorisations, but in triangular ones. Thus, we have
to modify induction base, which now becomes even easier, and
superficially the reduction step itself. It is a total mystery, why
we failed to notice these obvious modifications when writing
\cite{VSS}.
\par
What is truly amazing here, is that as in \cite{VSS} the usual
{\it linear\/} stable rank condition works for groups of all
types! Before, this decomposition was known for the
Chevalley group $G(\Phi,R)$ itself, under the following
[stronger!] assumptions on $R$.
\par\smallskip
$\bullet$ For $\Phi=\A_l,\C_l$ under $\sr(R)=1$.
\par\smallskip
$\bullet$ For $\Phi=\B_l,\D_l$ under $\asr(R)=1$, or under
an appropriate unitary/form ring stable rank condition
$\Lambda\!\sr(R)=1$, $\sr(R,\Lambda)=1$, etc.
\par\smallskip
$\bullet$ For exceptional groups, when $R$ is semi-local.
\par\smallskip
\noindent These results, especially for classical groups, were
immediately obvious, after the introduction of the corresponding
stability conditions by Hyman Bass \cite{bass} and by Michael Stein
\cite{stein73}\footnote{Actually, absolute stable rank as such was
first introduced by David Estes and Jack Ohm \cite{EO}, but its role
in the proof of stability results for orthogonal groups was noted
only by Stein.}. Under these conditions Gauss decomposition for
classical groups was [re]discovered dozens of times, and in the last
section we provide assorted references.
\par
Our Theorem 1 divorces existence of Gauss decomposition
from the triviality of $K_1(\Phi,R)$. In fact, it shows
that these stronger stability conditions are only needed
to ensure that $G_{\sic}(\Phi,R)=E_{\sic}(\Phi,R)$, but
are not necessary for the {\it elementary\/} Chevalley
group $E(\Phi,R)$ to admit Gauss decomposition!
\par
Let us state some immediate corollaries of Theorem~1.
\begin{proclaim}
{Corollary 1} Let\/ $\Phi$ be a reduced irreducible root system
and\/ $R$ be a commutative ring such that $\sr(R)=1$. Then any
element $g$ of the elementary Chevalley group\/ $E(\Phi,R)$
is conjugate to an element of
$$ U(\Phi,R)H(\Phi,R)U^-(\Phi,R). $$
\end{proclaim}
\begin{proclaim}
{Corollary 2} Let\/ $\Phi$ be a reduced irreducible root system
and\/ $R$ be a commutative ring such that $\sr(R)=1$. Then
the elementary Chevalley group\/ $E(\Phi,R)$ admits
unitriangular factorisation
$$ E(\Phi,R)=U(\Phi,R)U^-(\Phi,R)U(\Phi,R)U^-(\Phi,R)U(\Phi,R) $$
\noindent
of length\/ $5$.
\end{proclaim}
Notice that this corollary is a very broad generalisation
of results on unitriangular factorisations obtained by
Martin Liebeck, Laszlo Pyber, Laszlo Babai and Nikolay
Nikolov \cite{LP,BNP}, with a {\it terribly\/} much easier proof.
Actually, Theorem~1 of \cite{VSS}, which is proven by
essentially the same method, but starts with a slightly
more precise induction base, asserts that under condition
$\sr(R)=1$ the elementary Chevalley group $E(\Phi,R)$
admits unitriangular factorisation
$$ G(\Phi,R)=U(\Phi,R)U^-(\Phi,R)U(\Phi,R)U^-(\Phi,R) $$
\noindent
of length 4.
\par
The present note is a by-product of our joint work on arithmetic
problems of our cooperative Russian--Indian project ``Higher
composition laws, algebraic $\K$-theory and exceptional groups'' at
the Saint Petersburg State University, Tata Institute of Fundamental
Research (Mumbai) and Indian Statistical Institute (Bangalore).
\par
In \S~1 we recall some fundamentals concerning $\sr(R)=1$ and
similar stability conditions. In \S\S~2 and~3 we introduce basic
notation related to Chevalley groups and their parabolic subgroups.
In \S~4 we prove another version of Tavgen's rank reduction theorem,
which immediately implies Theorem 1. Finally, in \S~5 we discuss
existing literature on the subject and state several unsolved
problems.
\par
The proofs in the present paper, as also in \cite{VSS}, are based
on a slight variation of an idea by Oleg Tavgen. We finished this
paper at the end of July 2011, and just as we planned to send it
to Oleg, we were deeply shocked and grieved by the news of his
sudden and untimely death. We dedicate this paper to his memory.


\section{Stability conditions}

Recall that a ring $R$ has {\bf stable rank $1$\/}, if for
all $x,y\in R$, which generate $R$ as a {\it right\/} ideal,
there exists a $z\in R$ such that $x+yz$ is {\it right\/}
invertible. In this case we write~$\sr(R)=1$.
\par
It is classically known that rings of stable rank 1 are actually
{\it weakly finite\/} (Kaplansky---Lenstra theorem), so that in
their definition one could from the very start require
that $x+yz\in R^*$. Since for the linear case the result is
well known, and Chevalley groups of other types only exist over
commutative rings, from here on we assume that the ring $R$ is
commutative, in which case the proof below at the same time
demonstrates that $\SL(2,R)=E(2,R)$.

\par
Special linear groups are simply connected Chevalley groups of type
$\A_l$. In particular, $\SL(2,R)$ and $E(2,R)$ are the simply
connected Chevalley group of type $\A_1$ and its elementary
subgroup. All other notations are modified accordingly. Thus,
$$ U(2,R)=\begin{pmatrix} 1&*\\ 0&1\\ \end{pmatrix},
\qquad
U^-(2,R)=\begin{pmatrix} 1&0\\ *&1\\ \end{pmatrix} $$
\noindent
refers to the groups $U(\A_1,R)$ and $U^-(\A_1,R)$, under
the above identification of $G_{\sic}(\A_1,R)=\SL(2,R)$,
etc. Similarly, $H(\A_1,R)=T(\A_1,R)$ is now identified with
$$ H(2,R)=\left\{\begin{pmatrix} \e&0\\ 0&\e^{-1}\\ \end{pmatrix}
\bigg\vert\ \e\in R^*\right\}, $$
\noindent
in other words, with the group of diagonal matrices with
determinant 1, usually denoted by $\SD(2,R)$.
\par
The proof of the following lemma is essentially contained
already in \cite{bass} and was rediscovered {\it dozens\/}
of times after that. The lemma itself, induction base,
is the only step in the proof of Theorem~1 that invokes
stability condition.
\begin{proclaim}
{Lemma 1} Let\/ $R$ be a commutative ring of stable
rank\/ $1$. Then
$$ \SL(2,R)=E(2,R)=U(2,R)H(2,R)U^-(2,R)U(2,R). $$
\end{proclaim}
\begin{proof}[Proof]
Consider an arbitrary matrix
$g=\begin{pmatrix} a&b\\ c&d\\ \end{pmatrix}\in\SL(2,R)$.
Since rows of an invertible matrix are unimodular, one has
$cR+dR=R$. Since $\sr(R)=1$, there exists such an $z\in R$,
that $d+cz\in R^*$. Thus,
$$ \begin{pmatrix} 1&-(b+az)(d+cz)^{-1}\\ 0&1\\ \end{pmatrix}
\begin{pmatrix} a&b\\ c&d\\ \end{pmatrix}
\begin{pmatrix} 1&z\\ 0&1\\ \end{pmatrix}=
\begin{pmatrix} (d+cz)^{-1}&0\\ c&d+cz\\ \end{pmatrix}, $$
\noindent
as claimed.
\end{proof}


\section{Chevalley groups}

\def\wP{\mathcal P}
\def\wQ{\mathcal Q}

Our notation pertaining to Chevalley groups are utterly standard
and coincide with the ones used in \cite{vavilov91,VP96},
where one can find many further references.
\par
Let as above $\Phi$ be a reduced irreducible root system of rank
$l$, $W=W(\Phi)$ be its Weyl group and $\wP$ be a weight lattice
intermediate between the root lattice $\wQ(\Phi)$ and the weight
lattice $\wP(\Phi)$. Further, we fix an order on $\Phi$ and denote
by $\Pi=\{\a_1,\ldots,\a_l\}$, $\Phi^+$ and $\Phi^-$ the corresponding
sets of fundamental, positive and negative roots, respectively.
Our numbering of the fundamental roots follows Bourbaki.
Finally, let $R$ be a commutative ring with 1, as usual,
$R^*$ denotes its multiplicative group.
\par
It is classically known that with these data one can associate
the Chevalley group $G=G_{\wP}(\Phi,R)$, which is the group of
$R$-points of an affine groups scheme $G_{\wP}(\Phi,-)$, known as
the Chevalley---Demazure group scheme. In the case
$\wP=\wP(\Phi)$ the group $G$ is called simply connected and
is denoted by $G_{\sic}(\Phi,R)$. In the opposite case
$\wP=\wQ(\Phi)$ the group $G$ is called adjoint and
is denoted by $G_{\ad}(\Phi,R)$. Many results do not depend
on the lattice $\wP$ and hold for all groups of a given
type $\Phi$. In such cases we omit any reference to
$\wP$ in the notation and denote by $G(\Phi,R)$
{\it any} Chevalley group of type $\Phi$ over $R$.
\par
In what follows, we fix a split maximal torus $T(\Phi,-)$
of the group scheme $G(\Phi,-)$ and set $T=T(\Phi,R)$.
As usual, $X_{\a}$, $\a\in\Phi$, denotes a unipotent root subgroup
in $G$, elementary with respect to $T$. We fix isomorphisms
$x_{\a}:R\mapsto X_{\a}$, so that
$X_{\a}=\{x_{\a}(\xi)\mid\xi\in R\}$, which are interrelated
by the Chevalley commutator formula, see \cite{carter,steinberg}. Further,
$E(\Phi,R)$ denotes the elementary subgroup of $G(\Phi,R)$,
generated by all root subgroups $X_{\a}$, $\a\in\Phi$.
\par
Elements $x_{\a}(\xi)$, $\a\in\Phi$, $\xi\in R$, are called
[elementary] unipotent root elements or, for short, simply
root unipotents. Next, let $\a\in\Phi$ and  $\e\in R^{*}$.
As usual, we set
$w_{\a}(\e)=x_{\a}(\e)x_{-\a}(-\e ^{-1})x_{\a}(\e)$ and
$h_{\a}(\e)=w_{\a}(\e)w_{\a}(1)^{-1}$. Elements $h_{\a}(\e)$
are called semisimple root elements. Define
$$ H(\Phi,R)=\langle h_{\a}(\e),\ \a\in\Phi,\ \e\in R^*\rangle. $$
\noindent
For a simply connected group one has
$$ H_{\sic}(\Phi,R)=T_{\sic}(\Phi,R)=\Hom(\wP(\Phi),R^*). $$
\noindent
In general, though, $H(\Phi,R)=T(\Phi,R)\cap E(\Phi,R)$ can be
--- and even over a field usually is! --- strictly smaller,
than $T(\Phi,R)$.
\par
Finally, let $N=N(\Phi,R)$ be the algebraic normaliser of the
torus $T=T(\Phi,R)$, i.~e.\ the subgroup, generated by $T=T(\Phi,R)$
and all elements $w_{\a}(1)$, $\a\in\Phi$. The factor-group
$N/T$ is canonically isomorphic to the Weyl group $W$, and for each
$w\in W$ we fix its preimage $n_{w}\in N$. Clearly, such a
preimage can be taken in $E(\Phi,R)$. Indeed, for a root
reflection $w_{\a}$ one can take $w_{\a}(1)\in E(\Phi,R)$
as its preimage, any element $w$ of the Weyl group can
be expressed as a product of root reflections.
In particular, we get the following classical result.

\begin{proclaim}{Lemma 2}
The elementary Chevalley group\/ $E(\Phi,R)$ is generated by
unipotent root elements\/ $x_{\a}(\xi)$, $\a\in\pm\Pi$, $\xi\in R$,
corresponding to the fundamental and negative fundamental roots.
\end{proclaim}

Further, let $B=B(\Phi,R)$ and $B^-=B^-(\Phi,R)$ be a pair of
opposite Borel subgroups containing $T=T(\Phi,R)$, standard
with respect to the given order. Recall that $B$ and $B^-$
are semidirect products $B=T\rightthreetimes U$ and
$B^-=T\rightthreetimes U^-$, of the torus $T$ and their
unipotent radicals
$$
\begin{aligned}
U&=U(\Phi,R)=
\big\langle x_\a(\xi),\ \a\in\Phi^+,\ \xi\in R\big\rangle, \\
\noalign{\vskip 4pt}
U^-&=U^-(\Phi,R)=
\big\langle x_\a(\xi),\ \a\in\Phi^-,\ \xi\in R\big\rangle. \\
\end{aligned}
$$
\noindent
Here, as usual, for a subset $X$ of a group $G$ one denotes by
$\langle X\rangle$ the subgroup in $G$ generated by $X$.
Semidirect product decomposition of $B$ amounts to saying that
$B=TU=UT$, and at that $U\trianglelefteq B$ and $T\cap U=1$.
Similar facts hold with $B$ and $U$ replaced by $B^-$ and $U^-$.
Sometimes, to speak of both subgroups $U$ and $U^-$ simultaneously,
we denote $U=U(\Phi,R)$ by $U^+=U^+(\Phi,R)$.
\par
In general, one can associate a subgroup $E(S)=E(S,R)$
to any closed set $S\subseteq\Phi$. Recall that a subset
$S$ of $\Phi$ is called {\it closed\/}, if for any two roots
$\a,\b\in S$ the fact that $\a+\b\in\Phi$, implies that already
$\a+\b\in S$. Now, we define $E(S)=E(S,R)$ as the subgroup generated by
all elementary root unipotent subgroups $X_{\a}$, $\a\in S$:
$$ E(S,R)=\langle x_{\a}(\xi),\quad \a\in S,\quad \xi\in R\rangle. $$
\noindent In this notation, $U$ and $U^{-}$ coincide with
$E(\Phi^{+},R)$ and $E(\Phi^{-},R)$, respectively. The groups
$E(S,R)$ are particularly important in the case where $S$ is a {\it
special\/} (= {\it unipotent\/}) set of roots; in other words, where
$S\cap(-S)=\varnothing$. In this case $E(S,R)$ coincides with the
{\it product\/} of root subgroups $X_{\a}$, $\a\in S$, in some/any
fixed order.
\par
Let again $S\subseteq\Phi$ be a closed set of roots. Then $S$ can be
decomposed into a disjoint union of its {\it reductive\/} (= {\it
symmetric\/}) part $S^{r}$, consisting of those $\a\in S$, for which
$-\a\in S$, and its {\it unipotent\/} part $S^{u}$, consisting of
those $\a\in S$, for which $-\a\not\in S$. The set $S^{r}$ is a
closed root subsystem, whereas the set $S^{u}$ is special. Moreover,
$S^{u}$ is an {\it ideal\/} of $S$, in other words, if $\a\in S$,
$\b\in S^{u}$ and $\a+\b\in\Phi$, then $\a+\b\in S^{u}$. {\it Levi
decomposition\/} asserts that the group $E(S,R)$ decomposes into
semidirect product $E(S,R)=E(S^r,R)\rightthreetimes E(S^u,R)$ of its
{\it Levi subgroup\/} $E(S^{r},R)$ and its {\it unipotent
radical\/}~$E(S^{u},R)$.

\section{Elementary parabolic subgroups}

\noindent
The main role in the proof of Theorem 1 is played by Levi
decomposition for elementary parabolic subgroups. Denote by
$m_k(\a)$ the coefficient of $\a_k$ in the expansion of $\a$
with respect to the fundamental roots:
$$ \a=\sum m_k(\a)\a_k,\quad 1\le k\le l. $$
\par
Now, fix an $r=1,\ldots,l$ -- in fact, in the reduction to smaller
rank it suffices to employ only terminal parabolic subgroups,
even only the ones corresponding to the first and the last
fundamental roots, $r=1,l$. Denote by
$$ S=S_r=\big\{\a\in\Phi,\ m_r(\a)\geq 0\big\} $$
\noindent
the $r$-th standard parabolic subset in $\Phi$. As usual,
the reductive part $\Delta=\Delta_r$ and the special part
$\Sigma=\Sigma_r$ of the set $S=S_r$ are defined as
$$ \Delta=\big\{\alpha\in\Phi,\ m_r(\alpha) = 0\big\},\quad
\Sigma=\big\{\alpha\in\Phi,\ m_r(\alpha) > 0\big\}. $$
\noindent
The opposite parabolic subset and its special part are defined
similarly
$$ S^-=S^-_r=\big\{\alpha\in\Phi,\ m_r(\alpha)\leq 0\big\},\quad
\Sigma^-=\big\{\alpha\in\Phi,\ m_r(\alpha)<0\big\}. $$
\noindent
Obviously, the reductive part $S^-_r$ equals $\Delta$.
\par
Denote by $P_r$ the {\it elementary\/} maximal parabolic
subgroup of the elementary group $E(\Phi,R)$. By definition,
$$ P_r=E(S_r,R)=\big\langle x_\alpha(\xi),\ \alpha\in S_r,
\ \xi\in R \big\rangle. $$
\noindent
Now Levi decomposition asserts that the group $P_r$ can be represented
as the semidirect product
$$ P_r=L_r\rightthreetimes U_r=E(\Delta,R)\rightthreetimes E(\Sigma,R) $$
\noindent
of the elementary Levi subgroup $L_r=E(\Delta,R)$ and the unipotent
radical $U_r=E(\Sigma,R)$. Recall that
$$ L_r=E(\Delta,R)=\big\langle x_\alpha(\xi),\quad \alpha\in\Delta,
\quad \xi\in R \big\rangle, $$
\noindent
Whereas
$$ U_r=E(\Sigma,R)=
\big\langle x_\alpha(\xi),\ \alpha\in\Sigma,\ \xi\in R\big\rangle. $$
\noindent
A similar decomposition holds for the opposite parabolic subgroup
$P_r^-$, whereby the Levi subgroup is the same as for $P_r$,
but the unipotent radical $U_r$ is replaced by the opposite unipotent
radical $U_r^-=E(-\Sigma,R)$
\par
As a matter of fact, we use Levi decomposition in the following
form. It will be convenient to slightly change the notation
and write $U(\Sigma,R)=E(\Sigma,R)$ and $U^-(\Sigma,R)=E(-\Sigma,R)$.

\begin{proclaim}{Lemma 3}
The group\/ $\big\langle U^{\sigma}(\Delta,R),U^\rho(\Sigma,R)\big\rangle$,
where\/ $\sigma,\rho=\pm 1$, is the semidirect product of its
normal subgroup\/ $U^\rho(\Sigma,R)$ and the complementary subgroup\/
$U^{\sigma}(\Delta,R)$.
\end{proclaim}

In other words, it is asserted here that the subgroups
$U^{\pm}(\Delta,R)$ normalise each of the groups
$U^{\pm}(\Sigma,R)$, so that, in particular, one has the following
four equalities for products
$$ U^{\pm}(\Delta,R)U^{\pm}(\Sigma,R)=U^{\pm}(\Sigma,R)U^{\pm}(\Delta,R), $$
\noindent
and, furthermore, the following four obvious equalities for
intersections hold:
$$ U^{\pm}(\Delta,R)\cap U^{\pm}(\Sigma,R)=1. $$
\par
In particular, one has the following decompositions:
$$ U(\Phi,R)=U(\Delta,R)\rightthreetimes U(\Sigma,R),
\quad
U^-(\Phi,R)=U^-(\Delta,R)\rightthreetimes U^-(\Sigma,R). $$


\section{Proof of Theorem~1}

The following result, like Theorem~3 of \cite{VSS}, is another
minor elaboration of Proposition 1 from the paper by Oleg
Tavgen \cite{tavgen90}. Tavgen considered {\it unitriangular
factorisations\/}, in other words, expressions of
$E(\Phi,R)$ as products of $U(\Phi,R)$ and $U^-(\Phi,R)$,
$$ E(\Phi,R)=U(\Phi,R)U^-(\Phi,R)\ldots U^{\pm}(\Phi,R). $$
\noindent
Here, we are interested in {\it triangular factorisations\/},
in other words expressions of $E(\Phi,R)$ as products of
$$ B(\Phi,R)\cap E(\Phi,R)=H(\Phi,R)U(\Phi,R) $$
\noindent
and
$$ B^-(\Phi,R)\cap E(\Phi,R)=H(\Phi,R)U^-(\Phi,R). $$
\noindent
However, since $T(\Phi,R)$ --- and, a fortiori, $H(\Phi,R)$
--- normalises $U(\Phi,R)$ and $U^-(\Phi,R)$, we can collect
all toral factors together, and consider factorisations of
the form
$$ E(\Phi,R)=H(\Phi,R)U(\Phi,R)U^-(\Phi,R)\ldots
U^{\pm}(\Phi,R). $$
\noindent
The length of such a decomposition is the number of distinct
{\it triangular\/} factors, in other words, the number of
$U^{\pm}(\Phi,R)$ occuring in this product.

\begin{proclaim}{Theorem 2}
Let\/ $\Phi$ be a reduced irreducible root system of rank $l\ge 2$,
and\/ $R$ be a commutative ring. Suppose that for the two
subsystems\/ $\Delta=\Delta_1,\Delta_l$, the elementary Chevalley
group\/ $E(\Delta,R)$ admits a triangular factorisation
$$ E(\Delta,R)=H(\Delta,R)U(\Delta,R)U^-(\Delta,R)\ldots
U^{\pm}(\Delta,R) $$
\noindent
of length $L$. Then the elementary Chevalley
group\/ $E(\Phi,R)$ admits triangular factorisation
$$ E(\Phi,R)=H(\Phi,R)U(\Phi,R)U^-(\Phi,R)\ldots
U^{\pm}(\Phi,R) $$
\noindent
of the same length\/ $L$.
\end{proclaim}

The leading idea of Tavgen's proof is very general and
beautiful, and works in many other related situations.
It relies on the fact that for systems of rank $\ge 2$
every fundamental root falls into the subsystem of smaller
rank obtained by dropping either the first or the last
fundamental root. Similar consideration was used by
Eiichi Abe and Kazuo Suzuki \cite{abe} and \cite{AS} to extract
root unipotents in their description of normal subgroups.
Compare also the simplified proof of Gauss decomposition
with prescribed semisimple part by Vladimir Chernousov,
Erich Ellers, and Nikolai Gordeev \cite{CEG}.
\par\smallskip\noindent
{\bf Remark.} As was pointed out by the referee, the argument
below applies without any modification in a much more general
setting. Namely, it suffices to assume that the required
decomposition holds for elementary Chevalley groups
$E(\Delta,R)$ for {\it some\/} subsystems $\Delta\le\Phi$,
whose union contains all fundamental roots. These subsystems
do not have to be terminal, or even irreducible, for that
matter. However, we do not see any immediate application of
this more general form of Theorem 2, since the use of terminal
subsystems invariably gives better results depending on
weaker stability conditions.
\par\smallskip
Let us reproduce the details of the argument. By definition
$$ Y=H(\Phi,R)U(\Phi,R)U^-(\Phi,R)\ldots
U^{\pm}(\Phi,R) $$
\noindent
is a {\it subset\/} in $E(\Phi,R)$. Usually, the easiest way
to prove that a subset $Y\subseteq G$ coincides with the whole
group $G$ consists in the following.

\begin{proclaim}{Lemma 4}
Assume that\/ $Y\subseteq G$, $Y\neq \varnothing$,
and\/ $X\subseteq G$ be a symmetric generating set. If\/
$XY\subseteq Y$, then\/ $Y=G$.
 \end{proclaim}

\begin{proof}[Proof of Theorem $2$]
By Lemma 2 the group $G$ is generated by the fundamental
root elements
$$ X=\big\{x_{\a}(\xi)\mid \a\in\pm\Pi,\ \xi\in R\big\}. $$
\noindent
Thus, by Lemma 4 is suffices to prove that $XY\subseteq Y$.
\par
Let us fix a fundamental root unipotent $x_{\a}(\xi)$.
Since $\rk(\Phi)\ge 2$, the root $\a$ belongs to at least one
of the subsystems $\Delta=\Delta_r$, where $r=1$ or $r=l$,
generated by all fundamental roots, except for the first or the
last one, respectively. Set $\Sigma=\Sigma_r$ and express
$U^{\pm}(\Phi,R)$ in the form
$$ U(\Phi,R)=U(\Delta,R)U(\Sigma,R),\quad
U^-(\Phi,R)=U^-(\Delta,R)U^-(\Sigma,R). $$
\par
Using Lemma 3 we see that
\begin{multline*}
Y=H(\Phi,R)U(\Delta,R)U^-(\Delta,R)\ldots U^{\pm}(\Delta,R)\cdot\\
U(\Sigma,R)U^-(\Sigma,R)\ldots U^{\pm}(\Sigma,R).
\end{multline*}
\noindent
Since $\a\in\Delta$, one has $x_{\a}(\xi)\in E(\Delta,R)$, so that
the inclusion $x_{\a}(\xi)Y\subseteq Y$ immediately follows from
the assumption.
\end{proof}


\section{Final remarks}

A major application of Theorem~1 we have in mind, is the
commutator width of elementary Chevalley group.
\par
One of the major recent advances was the positive solution
of Ore's conjecture, asserting that any element of
$E_{\ad}(\Phi,K)$ over a field $K$ is a single commutator,
whenever this group is simple [as an abstract group]. For
large fields, say, all fields containing $\ge 8$ elements,
this was proven by Erich Ellers and Nikolai Gordeev \cite{EG},
using their remarkable results on Gauss decomposition with
prescribed semi-simple part, see \cite{CEG} and references there.
For small fields, this result was obtained by Martin Liebeck,
Eamond O'Brien, Aner Shalev and Pham Huu Tiep \cite{LOST},
using very delicate character estimates. It is essential that
the groups are {\it adjoint\/}. Beware, that in general
one may need {\it two\/} commutators to expess some
elements of $E_{\sic}(\Phi,K)$.
\par
We believe that solution of the following two problems is now
at hand. Compare the works of Arlinghaus, Leonid Vaserstein,
Ethel Wheland and You Hong \cite{VW1,VW2,you,AVY}, where this
is essentially done for classical groups, over rings subject
to $\sr(R)=1$ or some stronger stability conditions, and the
work by Nikolai Gordeev and Jan Saxl \cite{GS}, where this
is essentially done over local rings.
\begin{proclaim}
{Problem 1}
Under assumption\/ $\sr(R)=1$ prove that any element of\/
$E_{\ad}(\Phi,R)$ is a product of\/ $\le 2$ commutators in\/
$G_{\ad}(\Phi,R)$.
\end{proclaim}
\begin{proclaim}
{Problem 2}
Under assumption\/ $\sr(R)=1$ prove that any element of\/
$E(\Phi,R)$ is a product of\/ $\le 3$ commutators in\/
$E(\Phi,R)$.
\end{proclaim}
It may well be that under this assumption the commutator
width of $E(\Phi,R)$ is always $\le 2$, but so far we
were unable to control details concerning semisimple factors.
\par
It seems, that one can apply the same argument
to higher stable ranks. Solution of the following problem
would be a generalisation of \cite{DV88}, Theorem~4.
\begin{proclaim}
{Problem 3}
If the stable rank\/ $\sr(R)$ of\/ $R$ is finite, and for
some $m\ge 2$ the elementary linear group\/
$E(m,R)=E_{\sic}(A_{m-1},R)$ has bounded word length with
respect to elementary generators, then for all\/ $\Phi$
of sufficiently large rank one has
$$ E(\Phi,R)=\big(U(\Phi,R)U^-(\Phi,R)\big)^3. $$
\end{proclaim}
In particular, it would follow that in this case any element
of $E(\Phi,R)$ is a product of $\le 6$ commutators. In
fact, we expect a much better result.
\begin{proclaim}
{Problem 4}
If the stable rank\/ $\sr(R)$ of\/ $R$ is finite, and for
some $m\ge 2$ the elementary linear group\/
$E(m,R)$ has bounded word length with
respect to elementary generators, then for all\/ $\Phi$
of sufficiently large rank any element of\/
$E(\Phi,R)$ is a product of\/ $\le 4$ commutators in\/
$E(\Phi,R)$.
\end{proclaim}
Theorem~2 of \cite{VSS} is the first step towards construction
of short triangular factorisations of Chevalley groups
over Dedekind rings of arithetic type. At present, sharp
bounds depend on the Generalised Artin Conjecture, which in turn
depends on the Generalised Riemann Hypothesis, but with \cite{MP}
there is some hope to divorce these bounds from GRH and we
are presently working on that. This would
then be a crucial advance in the direction of the following result.
\begin{proclaim}
{Problem 5}
Let $R$ be a Dedekind ring of arithmetic type with infinite
multiplicative group. Prove that any element of\/
$E_{\ad}(\Phi,R)$ is a product of\/ $\le 3$ commutators in\/
$G_{\ad}(\Phi,R)$.
\end{proclaim}
Some of our colleagues expressed belief that any element
of $\SL(n,\Int)$, $n\ge 3$, is a product of $\le 2$ commutators.
However, for Dedekind rings with {\it finite\/} multiplicative
groups, such as $\Int$, at present we do not envisage any
{\it obvious\/} possibility to improve the generic bound $\le 4$
even for large values of $n$. Expressing elements of $\SL(n,\Int)$
as products of 2 commutators, if it can be done at all, should
require a lot of specific case by case analysis.
\par
Triangular factorisations, such as Gauss decomposition
considered in this paper, are the simplest instance of
{\it parabolic\/} factorisations. Recently, Sergei Sinchuk and
the third author obtained analogues of Dennis---Vaserstein
decomposition for arbitary pairs of maximal parabolic
subgroups $(P_r,P_s)$, $r<s$, in classical groups and pairs
of {\it terminal\/} parabolic subgroups in exceptional groups,
see \cite{SiV,VS10,VS11}. In some cases stronger conditions
than $\sr(R)<s-r$ were imposed. Now, it seems, that these
stronger conditions were only used to ensure surjective
stability of $K_1$ and are not needed to get decompositions
of the elementary group itself. Here, $U_{rs}^-=U_r^-\cap U_s^-$.
\begin{proclaim}
{Problem 6}
Prove Dennis---Vaserstein type decomposition
$$ G=P_rU_{rs}^-P_s $$
\noindent
for elementary Chevalley groups $E(\Phi,R)$, under restrictions
on the usual stable rank\/ $\sr(R)$.
\end{proclaim}
A similar problem for unitary groups was recently solved by
Sergei Sinchuk \cite{sinchuk}, as part of his efforts to improve
stability results for unitary $\K_1$, see \cite{BPT,BT}.
\par
Results of the present paper are also closely related to
Bass---Kolster type decompositions. The classical Bass---Kolster
decomposition for the group $\SL(n,R)$ has the form
$$ G=L_rU_rU^-_rU_rU^-_r=P_rU^-_rU_rU^-_r. $$
\noindent
For groups of other types, one has to vary one of the unipotent
radicals and gets decompositions of the form
$G=L_rU_rU^-_rU_sU^-_s$. In \cite{loos} Ottmar Loos addresses
the problem, whether one can shorten this decomposition. He
uses language of Jordan pairs, and his results apply to some
non-split reductive groups. In our context the problem he
studies amounts to asking, whether
$$ G=L_rU_rU^-_rU_r=P_rU^-_rU_r, $$
\noindent
for parabolic subgroups with Abelian unipotent radical.
He proves that for $\SL(n,R)$ existence of such a decomposition
is equivalent to $\sr(R)=1$, \cite{loos}, Corollary 3.9, but
this follows already from \cite{bass}. To establish similar
shorter Bass---Kolster type decompositions for groups of other
types, he imposes further conditions, similar in spirit
to von Neumann regularity. However, our Theorem 1 suggests that
no such conditions are necessary.
\begin{proclaim}
{Problem 7}
Prove Bass---Kolster type decompositions for elementary Chevalley
groups $E(\Phi,R)$, under restrictions on the usual stable
rank\/ $\sr(R)$.
\end{proclaim}
The authors would like to thank Alireza Abdollahi, Nikolai
Gordeev and You Hong for inspiring discussions of commutator
width and related problems, and Anastasia Stavrova, Alexei
Stepanov and Maxim Vsemirnov for some very pertinent remarks.
The third author thanks
Francesco Catino, Francesco de Giovanni and Carlo Scoppola for an
invitation to give a talk on commutator width and unitriangular
factorisations at the conference ``Advances in Group Theory
and Applications'' (Porto Cesareo -- June 2011), which helped him
to focus thoughts in this direction.


\end{document}